\documentclass[a4paper,11pt,reqno]{amsart}
\usepackage{graphicx}
\usepackage{amsmath}
\usepackage{amssymb}
\usepackage{oldgerm}
\usepackage{multicol}
\usepackage[all]{xy}
\newtheorem{te}{Theorem}[section]
\newtheorem{co}[te]{Corollary}

\newtheorem{de}[te]{\sc Definition}
\newtheorem{ex}[te]{Example}
\newtheorem{prop}[te]{Proposition}
\newtheorem{lm}[te]{Lemma}
\newtheorem{rem}[te]{Remark}

\begin{document}
\baselineskip=15pt
\title{Graded Blocks of Group Algebras With Dihedral Defect Groups}
\author{Dusko Bogdanic}
\address{Dusko Bogdanic \newline Mathematical Institute \\ University of Oxford \\  \newline 24-29 St.\
Giles \\ Oxford OX1 3LB \\ United Kingdom}

\email{bogdanic@maths.ox.ac.uk}
\date{\today}

\begin{abstract}
In this paper we investigate gradings on tame blocks of group algebras whose defect groups are dihedral. For this subfamily of tame blocks we classify gradings up to graded Morita equivalence, we transfer gradings via derived equivalences, and we check the existence, positivity and tightness of gradings. We classify gradings by computing the group of outer automorphisms that fix the isomorphism classes of simple modules. 
\end{abstract}
\maketitle \vspace{-10mm}
\section{Introduction}
\noindent In this paper we study gradings on tame blocks of group
algebras. Erdmann classified tame blocks of group algebras up to
Morita equivalence (cf.\ [\ref{ErK}]). A block of a group algebra
over a field of characteristic $p$ is of tame representation type
if and only if $p=2$ and its defect group is a dihedral,
semidihedral, or generalized quaternion group. If the defect group
of a block is a dihedral (respectively semidihedral, quaternion) group, then we say that the block is of dihedral (respectively semidihedral, quaternion) type. The number of simple
modules in a tame block is $1$, $2$ or $3$ (see [\ref{ErK}] for
more details).  Erdmann's
classification has been used by Holm to classify tame blocks up to
derived equivalence (the case of blocks with dihedral defect
groups and three simple modules has been dealt with by Linckelmann
in [\ref{Link}]). We will follow Erdmann's and Holm's classification,
and use some of the tilting complexes given in [\ref{H96}] and [\ref{Link}] to transfer gradings via derived equivalences in order to prove the existence of non-trivial gradings on an arbitrary dihedral block. 

As in the case of Brauer tree algebras (cf.\ [\ref{GBTA}]), we classify gradings up to graded Morita equivalence by computing the group of outer
automorphisms that fix the isomorphism classes of simple modules. From
our computation of these groups we are able to deduce that, in the
case of dihedral blocks with two simple modules, for different
scalars (which remain undetermined in Erdmann's classification) we
get algebras that are not derived equivalent.


The paper is organized as follows. In the second section we list some preliminary results that will be used throughout this paper. This section contains a classification criterion, and a criterion for tightness and positivity of gradings. In the third section we investigate gradings on dihedral blocks with three simple modules. In the fourth section we investigate gradings on dihedral blocks with two simple modules. The fifth section is devoted to dihedral blocks with one simple module. 

\subsection{Notation} Throughout this text $k$ will be an
algebraically closed field of  characteristic 2. All
algebras will be finite dimensional algebras over the field $k$, and all
modules will be left modules. The category of finite dimensional
$A$--modules is denoted by $A$--${\rm mod}$ and the full
subcategory of finite dimensional projective $A$--modules is
denoted by $P_A$. The derived category of bounded complexes over
$A$--${\rm mod}$ is denoted by $D^b(A)$, and the homotopy
category of bounded complexes over $P_A$ will be denoted by
$K^b(P_A)$.

\subsubsection{Graded modules} We say that an algebra $A$ is a graded algebra if $A$ is  the direct sum of subspaces $A=\bigoplus_{i\in\mathbb{Z}}
A_i$, such that $A_iA_j\subset A_{i+j}$, $i,j\in \mathbb{Z}$. If
$A_i=0$  for $i<  0$, we say that $A$ is positively graded. An
$A$-module $M$ is graded if it is the direct sum of subspaces
$M=\bigoplus_{i\in\mathbb{Z}} M_i$,  such that  $A_iM_j\subset
M_{i+j}$, for all $i,j\in \mathbb{Z}$. If $M$ is a graded
$A$--module, then $N=M\langle i\rangle$ denotes the graded module
given by $N_j=M_{i+j}$, $j\in \mathbb{Z}$. An $A$-module
homomorphism $f$ between two graded modules $M$ and $N$ is a
homomorphism of graded modules if $f(M_i)\subseteq N_i$, for all
$i\in \mathbb{Z}$. For a graded algebra $A$, we denote by
$A$--${\rm modgr}$ the category of graded finite dimensional
$A$--modules. We set ${\rm Homgr}_A(M,N):=\bigoplus_{i\in
\mathbb{Z}}{\rm Hom}_{A-gr}(M,N\langle i\rangle),$ where ${\rm
Hom}_{A-gr}(M,N\langle i \rangle)$ denotes the space of all graded
homomorphisms between $M$ and $N\langle i\rangle$ (the space of
homogeneous morphisms of degree $i$). There is an isomorphism of
vector spaces ${\rm Hom}_A(M,N)\cong {\rm Homgr}_A(M,N)$ that
gives us a grading on ${\rm Hom}_A(M,N)$ (cf.\ [\ref{NFV}],
Corollary 2.4.4.).


\subsubsection{Graded complexes} Let $X=(X^i,d^i)$ be a complex of $A$--modules. We say that $X$ is a complex of graded $A$--modules, or just a graded complex, if for
each $i\in \mathbb{Z}$, $X^i$ is a graded module and $d^i$ is a
homomorphism between graded $A$--modules. If $X$ is a graded
complex, then $X\langle j\rangle$ denotes the complex of graded
$A$--modules given by $(X\langle j\rangle)^i:=X^i\langle j\rangle$
and $d_{X\langle j\rangle}^i:=d^i$. Let $X$ and $Y$ be graded
complexes. A homomorphism $f=\{f^i\}_{i\in\mathbb{Z}}$ between
complexes $X$ and $Y$ is a homomorphism of graded complexes if for
each  $i\in \mathbb{Z}$, $f^i$ is a homomorphism of graded
modules. The category of complexes of graded $A$--modules will be
denoted by $C_{gr}(A)$. We set  ${\rm Homgr}_{A}(X, Y) :=
\bigoplus_{i\in\mathbb{Z}}{\rm Hom}_{C_{gr}(A)}(X,Y\langle
i\rangle ),$ where ${\rm Hom}_{C_{gr}(A)}(X,Y\langle i\rangle )$
denotes the space of graded homomorphisms between $X$ and
$Y\langle i\rangle$ (the space of homogeneous morphisms of degree
$i$). As for modules, we have an isomorphism of vector spaces
${\rm Homgr}_{A}(X, Y)\cong {\rm Hom}_{A}(X, Y)$ that gives us a
grading on ${\rm Hom}_{A}(X, Y)$. From this we get a grading on
${\rm Hom}_{K^b(A-{\rm mod})}(X, Y)$, since the subspace of zero
homotopic maps is homogeneous. We denote this graded space by
${\rm Homgr}_{K^b(A-{\rm mod})}(X, Y)$.

Unless otherwise stated, for a graded algebra $A$ given by a quiver and relations, we will assume
that the projective indecomposable $A$-modules are graded as in Example \ref{ex:quiverGrading} below, i.e.\ we will assume that their tops are in degree 0. 
We note here that if we have two different gradings on an indecomposable module (bounded complex), then they differ only by a shift (cf.\ [\ref{BGS}], Lemma 2.5.3).

\section{Preliminaries}
\subsection{Derived equivalences}
We say that two symmetric algebras $A$ and $B$ are derived
equivalent if their derived categories of bounded complexes are
equivalent. From Rickard's theory we know that $A$ and $B$ are
derived equivalent if and only if there exists a tilting complex
$T$ of projective $A$--modules such that ${\rm
End}_{K^b(P_A)}(T)\cong B^{op}$.  For more details on derived
categories and derived equivalences we recommend [\ref{KNG}]. 

We remind the reader that derived equivalent algebras share many common properties. Among these is the identity component ${\rm Out}^0(A)$ of the group of outer automorphisms (cf.\ [\ref{Saorin}], Theorem 17 or [\ref{Rou}], Theorem 4.6).

\subsection{Algebraic groups and a classification criterion}

For a finite dimensional $k$-algebra $A$, there  is a
correspondence between gradings on $A$ and homomorphisms of
algebraic groups from $\textbf{G}_m$ to ${\rm Aut}(A)$, where
$\textbf{G}_m$ is the multiplicative group $k^*$ of the field $k$.
For each grading $A=\bigoplus_{i\in \mathbb{Z}}A_i$ there is a
homomorphism of algebraic groups $\pi \, : \, \textbf{G}_m
\rightarrow {\rm Aut}(A)$ where an element $x\in k^*$ acts on
$A_i$ by multiplication by $x^i$ (see [\ref{Rou}], Section 5). If
$A$ is graded and $\pi$ is the corresponding homomorphism, we will
write $(A,\pi)$ to denote that $A$ is graded with grading $\pi$.

\begin{de}\label{grMorDef} Let $(A,\pi)$ and $(A,\pi^{\prime})$ be two gradings
on a finite dimensional $k$-algebra $A$,  and let
$S_1,S_2,\dots,S_r$ be the isomorphism classes of simple
$A$-modules. We say that $(A,\pi)$ and $(A,\pi^{\prime})$ are
graded Morita equivalent if there exist integers $d_{ij}$, where
$1\leq j\leq {\rm dim }\, S_i$ and $1\leq i\leq r$, such that the
graded algebras $(A,\pi^{\prime})$ and ${\rm
Endgr}_{(A,\pi)}(\bigoplus_{i,j}P_i\langle d_{ij}\rangle)^{\rm
op}$ are isomorphic, where $P_i$ denotes the projective cover of
$S_i$.
\end{de}
Note that two graded algebras are graded Morita equivalent if and
only if their categories of graded modules are equivalent.

Let $A=\bigoplus_{i\in \mathbb{Z}}A_i$ be a grading on $A$. If $r\in\mathbb{Z}$, then $A=\bigoplus_{i\in \mathbb{Z}}B_i$, where $B_{ri}:=A_i$, $i\in \mathbb{Z}$, and $B_i:=0$ for $r\nmid i$, is a grading on $A$. This procedure of multiplying (or dividing) each degree by the same integer is called rescaling.

We now give some background on algebraic groups (more details can be found in [\ref{Borel}]). An algebraic
torus is a linear algebraic group isomorphic to
$\textbf{G}_m^n=\textbf{G}_m\times\dots\times \textbf{G}_m$ ($n$
factors) for some $n\geq 1$. A maximal torus in an algebraic group
$G$ is a closed subgroup of $G$ which is a torus but is not
contained in any larger torus. Tori are contained in $G^0$, the
connected component of $G$ that contains the identity element. For
a given torus $T$,  a cocharacter of $T$ is a  homomorphism of
algebraic groups from $\textbf{G}_m$ to $T$. A cocharacter of an
algebraic group $G$ is a homomorphism of algebraic groups from
$\textbf{G}_m$ to $T$, where $T$ is a maximal torus of $G$. We say
that cocharacters $\pi$ and $\pi^{\prime}$ of $G$ are conjugate if
there exists $g \in G$ such that $\pi^{\prime}(x)=g\pi(x)g^{-1}$
for all $x\in \textbf{G}_m$. We see that a grading on a finite
dimensional algebra $A$ can be seen as a cocharacter $\pi \, : \,
\textbf{G}_m\rightarrow \, {\rm Aut}(A)$. We will use the same
letter $\pi$ to denote the corresponding cocharacter of ${\rm
Out}(A)$, which is given by composition of $\pi$ and the canonical
surjection.

The following proposition tells us how to classify all gradings on
$A$ up to graded Morita equivalence.

\begin{prop}[[\ref{Rou}{]}, Corollary 5.9]\label{cokarak}
Two basic graded algebras $(A,\pi)$ and $(A,\pi^{\prime})$ are graded
Morita equivalent if and only if the corresponding cocharacters
$\,\pi \, : \, \textbf{G}_m \rightarrow {\rm Out}(A)$  and
$\,\pi^{\prime} \, : \, \textbf{G}_m \rightarrow {\rm Out}(A)$ are
conjugate.
\end{prop}
From this proposition we see that in order to classify gradings on
$A$ up to graded Morita equivalence, we need to compute maximal
tori in  ${\rm Out}(A)$. Let ${\rm Out}^K (A)$ be  the subgroup of
${\rm Out}\, (A)$ of those automorphisms fixing the isomorphism
classes of simple $A$-modules. Since ${\rm Out}^K (A)$ contains
${\rm Out}^0(A)$, the connected component of ${\rm Out}(A)$ that
contains the identity element, we have that maximal tori in ${\rm
Out}(A)$ are actually contained in ${\rm Out}^K (A)$. It follows that it is sufficient to compute maximal tori in ${\rm Out}^K (A)$. 

\begin{lm}\label{KlGr}
Let $A$ be a basic finite dimensional algebra such that the
maximal tori in ${\rm Out}(A)$ are isomorphic to ${\textbf{G}}_m$.
Up to graded Morita equivalence and rescaling there is a unique
grading on $A$.
\end{lm}
\noindent{\bf Proof.} We saw at the beginning of this section that
gradings on $A$ correspond to cocharacters of ${\rm Aut}(A)$. If
$A=\bigoplus_{i\in \mathbb{Z}}A_i$ is a grading on $A$, then the
corresponding cocharacter is given by the action of $x$ on $A_i$
by $x*a_i=x^ia_i$, where $a_i\in A_i$.
Let $T$ and $T^{\prime}$ be two maximal tori in ${\rm Out}(A)$.
Let $\tau$ be a cocharacter of ${\rm Out}(A)$ such that its image
is contained in $T^{\prime}$. Since any two maximal tori in ${\rm
Out}(A)$ are conjugate, there exists an invertible element $a$
such that $aT^{\prime}a^{-1}=T$. The cocharacter given by
$x\mapsto a\tau(x)a^{-1}$, $x\in \mathbf{G}_m$, is conjugate to
$\tau$ and its image is contained in $T$. This cocharacter gives
rise to a grading which is graded Morita equivalent to the grading
given by $\tau$. It follows that when classifying gradings on $A$
up to graded Morita equivalence it is  sufficient to consider
cocharacters whose image is in $T$. The only homomorphisms from ${
\textbf{G}}_m$ to ${\textbf{G}}_m\cong T$ are given by maps
$x\mapsto x^r$, for $x\in {\textbf{G}}_m$ and $r\in \mathbb{Z}$.
Let $\pi\, :\, \textbf{G}_m\rightarrow {\rm Out}(A)$, $x\mapsto
x^l$, be the cocharacter that corresponds to the grading
$A=\bigoplus_{i\in \mathbb{Z}}A_i$. If we rescale this grading by
multiplying by $r\in \mathbb{Z}$, then we get the grading
$A=\bigoplus_{i\in \mathbb{Z}}B_i$, where $B_{ri}:=A_i$, $i\in
\mathbb{Z}$, and $B_i:=0$, for $r\nmid i$. This grading
corresponds to the cocharacter $\pi_1\, :\,
\textbf{G}_m\rightarrow {\rm Out}(A)$, $x\mapsto x^{rl}$. This is
easily seen if one thinks of the action of $x\in \mathbf{G}_m$ on
$B_{ri}$. If $b_{ri}\in B_{ri}$, then $b_{ri}=a_i$, $a_i\in A_i$.
The action of $x$ is given by
$$\pi_1(x)(b_{ri})=x*b_{ri}=x^{ri}b_{ri}=x^{ri}a_i=(\pi(x))^r(a_i).$$
We see that the grading corresponding to the cocharacter $x\mapsto
x^r$, $r\in \mathbb{Z}$, can be obtained by rescaling by $r$ from
the grading corresponding to the cocharacter $x\mapsto x$. It
follows that there is a unique grading on $A_{\Gamma}$ up to
rescaling (dividing or multiplying each degree by the same
integer) and graded Morita equivalence (shifting each projective
indecomposable module by an integer). $\blacksquare$

\subsection{A criterion for tightness and positivity}
\begin{prop}Let $A=\bigoplus_{i\geq0}A_i$ be
a positively graded algebra. Let $e$ and $f$ be homogeneous
primitive idempotents such that $Ae\cong Af$. Then $Ae$ and $Af$
are isomorphic as graded $A$-modules.
\end{prop}
\noindent {\bf Proof.} The modules $Ae=\bigoplus_{i\geq 0}A_ie$
and $Af=\bigoplus_{i\geq 0}A_if$ are positively graded. Since $Ae\cong Af$,
there exists an invertible element $a$ such that $aea^{-1}=f.$ If
$a_0$ is the degree 0 component of $a$, then $a_0ea^{-1}_0 = f$.
Right multiplication by $a_0$ is an isomorphism between the graded
modules $Af$ and $Ae$. $\blacksquare$

\begin{ex}\label{exProjGrad}
{\rm Let $A$ be a positively graded algebra and let $P$ be a projective
indecomposable $A$-module. There is a canonical way to grade $P$
as follows. \noindent Let $\{f_1,f_2,\dots,f_r\}$ be a complete
set of primitive orthogonal idempotents. If $e_i$ is the degree 0
component of $f_i$, then by comparing degree 0 components of
$f_i^2=f_i$, we conclude that $e_i$ is a primitive idempotent.
Hence,  $\{e_1,e_2,\dots, e_r\}$ is a complete set of primitive
orthogonal idempotents and $A=\bigoplus_{i=1}^r Ae_i$ is a sum of
graded modules. The projective indecomposable module  $P$ is
isomorphic to $Ae_i$ for some $i$. This gives us a grading on $P$,
which by the previous proposition does not depend on the choice of
the idempotent $e_i$. It follows that every projective $A$-module
is graded as a direct sum of graded modules.}
\end{ex}

\begin{de}
Let $A$ be a graded algebra. An ideal $I$ of $A$ is called
homogeneous if it is generated by homogeneous elements. 
\end{de}

\begin{lm}
Let $A$ be a graded algebra and let $I$ be a homogeneous ideal of
$A$. Then $A/I$ is a graded algebra.
\end{lm}
\noindent {Proof.} We define $(A/I)_i:=(A_i+I)/I$. $\blacksquare$
\begin{ex}\label{ex:quiverGrading}
{\rm Let $A$ be a finite dimensional algebra given by the quiver
$Q$ and the ideal of relations $I$, i.e.\ $A= kQ/I$.  The algebra
$kQ$ is generated, as an algebra, by the vertices and arrows of
$Q$. In order to grade $kQ$ it is sufficient to define the degrees
of the arrows since the vertices of $Q$ will be in degree 0. In
order to grade $kQ/I$, it is sufficient to ensure that $I$ is a
homogeneous ideal of $kQ$. In other words, if ${\rm
deg}(\alpha)={\rm deg}(\beta)$ for each relation $\alpha=\beta$
from  a generating set of $I$, where $\alpha$ and $\beta$ are
paths in $Q$ with the same source and the same target, then $I$ is
generated by homogeneous elements.

Let us assume that $A=kQ/I$ is graded in such a way that the
arrows and the vertices of $Q$ are homogeneous, and that $I$ is a
homogeneous ideal of $kQ$. Let $Ae$ be the projective indecomposable
module that corresponds to a vertex $e$ of the quiver $Q$. Then
$Ae$ is graded in a natural way as follows. As a vector space
$$Ae=\bigoplus_{\alpha}k\alpha,$$ where the sum  runs over the
non-zero paths $\alpha$ in the quiver $Q$ that have $e$ as their
target. If $\alpha$ is a path of length $l$, then the degree of
the $1$-dimensional subspace corresponding to $\alpha$ is $l$.
In this way, we can grade the projective indecomposable
$A$-modules even if the grading on $A$ is not positive. }
\end{ex}

\begin{de}[{[}\ref{CPS}{]}, {\rm Section 4}]\label{de:tight}
Let ${\rm gr}_{{\rm rad}\,A }(A)$ be the graded algebra given by
the radical filtration on a $k$-algebra $A$. We say that $A$ is a
tightly graded algebra if there is an algebra isomorphism
$$A\cong {\rm gr}_{{\rm rad}\,A
}(A).$$
\end{de}

By Proposition 4.4 in [\ref{CPS}], $A$ is tightly graded if and only if there exists a positive grading $A=\bigoplus_{i\geq 0} A_i$ such that $A_0$ is semisimple, and $A$ is generated, as an algebra, by $A_0$ and $A_1$. Such a grading is called tight.

\begin{lm} \label{lm: conjgrad}
Let $A=\bigoplus_{i\geq 0}A_i$ be a tight grading on a  $k$-algebra $A$.
If $a$ is an invertible element in $A$, then
$$A=\bigoplus_{i\geq 0}aA_ia^{-1}$$ is a tight grading on $A$.
\end{lm}
\noindent {\bf Proof.} This is obvious. $\blacksquare$

\begin{lm}
Let $A=\bigoplus_{i\geq 0}A_i$ be a tight grading on $A$.
If $A_{\geq i}:=\bigoplus_{j\geq i} A_j$, then
$${\rm rad}^i\,A=A_{\geq i},$$
and $A_0$ is a maximal semisimple subalgebra of $A$.
\end{lm}
\noindent {\bf Proof.} Since $A$ is an artinian algebra, $A_{\geq
1}$ is a nilpotent ideal. Hence, $A_{\geq 1}\subset {\rm rad}\, A$.
Let $S$ be a maximal semisimple subalgebra of $A$ such that $A=S\oplus
{\rm rad}\, A$.  Any two maximal
semisimple subalgebras of $A$ are conjugate (cf.\ [\ref{DrKir}], Theorem 6.2.1), and hence have the same dimension.
Because $A_0$ is a semisimple subalgebra, the dimension argument
gives us that $A_0$ is a maximal semisimple subalgebra and that
$A_{\geq 1}={\rm rad}\, A$. It follows easily that $A_{\geq
i}={\rm rad}^i\, A$, for $i\geq 1.\, \blacksquare$

\begin{lm}\label{lm: tight}
Let $A$ be an algebra given by the quiver $Q$ and the ideal of
relations $I$. Let $Q$ be  such that there are no
multiple arrows having the same source and the same target. If $A$
is a tightly graded algebra, then there exists a tight grading on
$A$ such that  for every arrow $\alpha$ of the quiver $Q$, there
exists a degree $1$ element $t_{\alpha}$ of the form
$\alpha+y_{\alpha}$, where $y_{\alpha}\in {\rm rad}^2\, A$ is a
linear combination of paths that have the same source and the same
target as $\alpha$.
\end{lm}
\noindent {\bf Proof.} Let us assume that
$A=\bigoplus_{i\geq 0}A_i$ is a tight grading on $A$. From
the previous lemma it follows that $A_0$ is a maximal semisimple  subalgebra
of $A$. Since any two maximal semisimple subalgebras are conjugate (cf.
[\ref{DrKir}], Theorem 6.2.1), by Lemma \ref{lm: conjgrad}, we can
assume that $A_0=S$, where $S$ is the maximal semisimple
subalgebra given by the linear span of the vertices of $Q$. Let
$\alpha$ be an arrow of $Q$ and let $x_1,\dots,x_s$ be degree 1
elements such that $\{x_1+{\rm rad}^2\, A,\dots,x_s+{\rm rad}^2\,
A\}$ is a basis of ${\rm rad}\, A/{\rm rad}^2\, A$. Then $\alpha$
can be written as a linear combination of homogeneous elements
$$\alpha=\sum_{i}\lambda_ix_i+y,$$ where $y\in {\rm rad}^2\, A.$
Because vertices are homogeneous, we can multiply this equation
from the left by $e_s$, the source vertex of $\alpha$, and from
the right by $e_t$, the target vertex of $\alpha$. We still get
$\alpha$ as a linear combination of homogeneous elements
$$\alpha=\sum_i\lambda_ie_sx_ie_t+e_sye_t.$$
By our assumption, the quiver $Q$ does not contain multiple arrows
with the same source and the same target. It follows that
$$\alpha=\sum_i\lambda_i(\mu_i \alpha+e_sz_ie_t)+e_sye_t,$$
where we assume that $x_i=\mu_i\alpha+w_i+z_i$, where $w_i$ is a
linear combination of arrows of $Q$ that are different from
$\alpha$, and $z_i$ is a linear combination of paths of length
greater than 1. It follows that $\sum_i\lambda_i\mu_i=1$ and that
the element $t_{\alpha}:=\alpha+\sum_i\lambda_i e_sz_ie_t$ is a
degree 1 element in $A$. $\blacksquare$

\begin{rem}\label{remark tight}{\rm The previous lemma can be
used to prove that certain algebras are not tightly graded, as in
the following case.

We keep the notation of the previous lemma. Let us assume that
one of the generators of $I$, say $v$,  is a linear combination of
paths such that at least two of them are of a different length.
We can assume that $v=\sum_{i=1}^r\lambda_ip_i$, where $p_i$ is a path of length $s_i$, i.e.\ $p_i=\alpha_{i1}\alpha_{i2}\cdots
\alpha_{is_i}$, where
$\alpha_{i1},\alpha_{i2},\dots,\alpha_{is_i}$ are arrows of $Q$.
If from the structure of $A$ it follows that
$p_i=t_{\alpha_{i1}}t_{\alpha_{i2}}\cdots t_{\alpha_{is_i}},$ for
all $p_i$,  then
${\rm deg}(t_{\alpha_{i1}}t_{\alpha_{i2}}\cdots
t_{\alpha_{is_i}})={\rm deg}(p_i)=s_i$, where $t_{\alpha}$ is as in the previous lemma. Without loss of
generality, let us assume that $p_{1}, p_{2},\dots,p_{m}$ are
paths of degree $s$, and that $p_{m+1}, \dots, p_r$ are
paths whose degree is greater than $s$. Then
$$\sum_{i=1}^m\lambda_ip_i=-\sum_{j=m+1}^r\lambda_jp_j.$$
Since the left-hand side of the above equality is a homogeneous
element of degree $s$, and the right-hand side is a sum of
homogeneous elements of degrees greater than $s$, we have a
contradiction, i.e.\ $A$ is not tightly graded.}
\end{rem}

Similar arguments can be used to prove that certain algebras given
by quivers and relations are not positively graded.
\begin{lm}\label{lm: tight2}
Let $A$ be an algebra given by the quiver $Q$ and the ideal of
relations $I$. Let $Q$ be such that there are no
multiple arrows having the same source and the same target. If $A$
is a positively graded algebra, then there exists a positive
grading on $A$ such that  for every arrow $\alpha$ of the quiver
$Q$, there exists a homogeneous element $t_{\alpha}$ of the form
$\alpha+y_{\alpha}$, where $y_{\alpha}\in {\rm rad}^2\, A$ is a
linear combination of paths that have the same source and the same
target as $\alpha$.
\end{lm}
\noindent {\bf Proof.} The arguments used in Example
\ref{exProjGrad} allow us to assume that the vertices of $Q$ are
homogeneous of degree 0. From the proof of  the previous
lemma it follows that for every arrow $\alpha$ there is a
homogeneous element $t_{\alpha}$ of the form $\alpha+y_\alpha$
such that its degree is non-negative. $\blacksquare$

\subsection{The group ${\rm Aut}(k[x]/(x^r))$.} 
This group will play an important role in our classification of gradings on dihedral blocks. We will denote it by $H_r$. 
\begin{de}\label{de: grupa}
We define $H_r$ to be the group $(k^*\times \underbrace{k\times k \times\dots\times
k}_{r-1}\, , \, *\, )$, where the multiplication $*$ is given by
\begin{equation} \label{eqGrupa}
\beta * \alpha:=\left(\sum_{i=1}^l\alpha_i\left(\sum_{\tiny \begin{array}{c}k_1+\dots+k_i=l\\
k_1,\dots,k_i>0\end{array}}\beta_{k_1}\beta_{k_2}\cdots\beta_{k_i}\right)\right)_{l=1}^r
\end{equation} 
\end{de}

Let $L$ be the subgroup of $H_r$ consisting of the elements of the
form $(1,\alpha_2,\dots, \alpha_r)$ and let $K$ be the subgroup of
$H_r$ consisting of the elements of the form
$(\alpha_1,0,\dots,0)$.

\begin{prop}
The group $H_r$ is a semidirect product of $L$ and $K$, where
$L\unlhd G$ is unipotent and the subgroup $K\cong
{\rm\textbf{G}}_m$ is a maximal torus in $H_r$.
\end{prop}
\noindent{\bf Proof.} This is straightforward. $\blacksquare$

\section{Three Simple Modules}

Any block with a dihedral defect group and three isomorphism
classes of simple modules is Morita equivalent to some algebra
from the following list (cf.\ [\ref{ErK}] or [\ref{H96}]).
\begin{itemize}
 \item [(1)]
For any $r\geq 1$, let $A_r$  be the algebra defined by
the quiver and relations
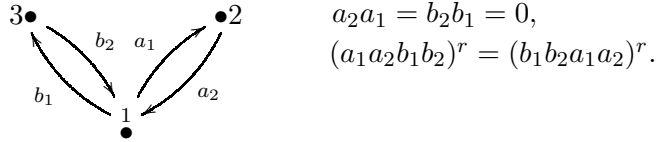
\begin{multicols}{2}
$
\xymatrix{&3\bullet\ar@/^/[dr]^{b_2}&&\bullet 2\ar@/^/[dl]^{a_2}\\
&&\stackrel{1}{\bullet}\ar@/^/[ur]^{a_1}\ar@/^/[ul]^{b_1}&}
$

$\!\!\!\!\!\!a_2a_1=b_2b_1=0,\\
(a_1a_2b_1b_2)^r=(b_1b_2a_1a_2)^r.$
\end{multicols}
\item [(2)] For any $r\geq 1$, let $B_r$ be the algebra defined by
the quiver and relations
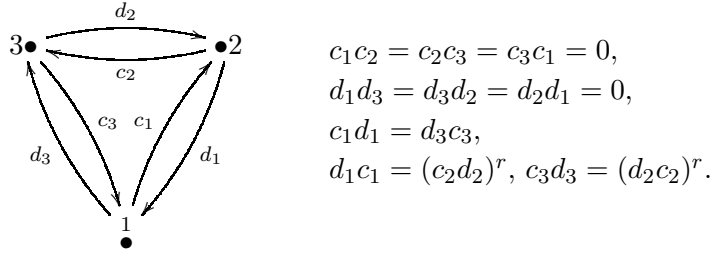
\begin{multicols}{2}
$
\xymatrix{&3\bullet\ar@/^/[ddr]^{c_3}\ar@/^/[rr]^{d_2}&&\bullet 2\ar@/^/[ll]^{c_2}\ar@/^/[ddl]^{d_1}\\&&&\\
&&\stackrel{1}{\bullet}\ar@/^/[uur]^{c_1}\ar@/^/[uul]^{d_3}&}
$
$ $\\ $ $\\ $ $\\ 
$\!\!\!\!\!\!c_1c_2=c_2c_3=c_3c_1=0,\\
d_1d_3=d_3d_2=d_2d_1=0,\\  
c_1d_1=d_3c_3, \\
d_1c_1=(c_2d_2)^r, \, c_3d_3=(d_2c_2)^r.\\
$
\end{multicols}
\item [(3)] For any $r\geq 2$, let $C_r$ be the algebra defined by the  quiver and relations
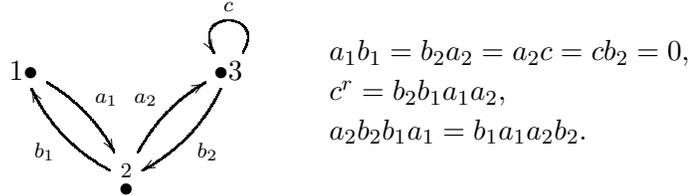
\begin{multicols}{2}
$
\xymatrix{&1\bullet\ar@/^/[dr]^{a_1}&&\bullet 3\ar@/^/[dl]^{b_2}\ar@(ur,ul)_{c}\\
&&\stackrel{2}{\bullet}\ar@/^/[ur]^{a_2}\ar@/^/[ul]^{b_1}&}
$

$ $\\ $\!\!\!\!\!\!\!a_1b_1=b_2a_2=a_2c=cb_2=0,$\\ $c^r=b_2b_1a_1a_2$,\\
$a_2b_2b_1a_1=b_1a_1a_2b_2.$
\end{multicols}
For $r=1$ we set $C_1=A_1$.
\end{itemize}

\subsection{Classification of gradings} We  start by classifying
all gradings up to graded Morita equivalence on $A_r$, $B_r$ and
$C_r$. In order to do this we need to compute  maximal tori in
${\rm Out}(A)$, where $A$ is $A_r$, $B_r$ or $C_r$. Since ${\rm
Out}^K(A)$, the group of outer isomorphisms that fix the
isomorphism classes of simple modules, contains ${\rm Out}^0(A)$,
and since ${\rm Out}^0(A)$ is invariant under derived equivalence
(cf.\ [\ref{Rou}], Theorem 4.6 or [\ref{Saorin}], Theorem 17), it
is sufficient to compute ${\rm Out}^K(A)$ for one of these
algebras. We will compute ${\rm Out}^K(C_r)$. Moreover, we will
see that  ${\rm Out}^K(C_r)$ and ${\rm Out}^0(C_r)$ are equal,
because ${\rm Out}^K(C_r)$ will turn out to be connected.

Let $\varphi$ be an arbitrary automorphism of $C_r$ fixing the isomorphism
classes of simple $C_r$-modules. The set $\{e_1, e_2, e_3\}$ of the vertices of the quiver of $C_r$ is a complete set of primitive orthogonal idempotents. Also, the set  $\{\varphi(e_1),\varphi(e_2),\varphi(e_3)\}$ is a complete set of  primitive orthogonal idempotents. From classical
ring theory (cf.\ [\ref{jac}], Theorem 3.10.2) we know that there exists an invertible element $x$ such that
$x^{-1}\varphi(e_i)x=e_{\sigma(i)}$, for all $i$, where $\sigma$
is some permutation. Since $\varphi$ fixes the isomorphism classes of
simple modules we can assume that $$\varphi(e_i)=e_i,\quad  i=1,2,3.$$ 

Since $\varphi({\rm rad}\,
C_r)\subseteq {\rm rad}\, C_r $, for a given arrow $t$ in the
quiver of $C_r$, $\varphi(t)$ is a linear combination of paths
whose source is the source of $t$ and whose target is the target
of $t$. It follows that
$$\begin{aligned}
  \varphi(a_1) &= \alpha_1a_1+\beta_1a_1a_2b_2, \\
  \varphi(a_2) &= \alpha_2a_2+\beta_2b_1a_1a_2,\\
  \varphi(b_1) &= \alpha_3b_1+\beta_3a_2b_2b_1, \\
  \varphi(b_2) &= \alpha_4b_2+\beta_4b_2b_1a_1, \\
  \varphi(c)&= \sum_{i=1}^r\gamma_ic^i,
\end{aligned}$$
where the $\alpha$'s, $\beta$'s and $\gamma$'s  are scalars. From
$a_1b_1=0$ and $b_2a_2=0$ we conclude that
$\alpha_1\beta_3+\alpha_3\beta_1=0$ and
$\alpha_4\beta_2+\alpha_2\beta_4=0$. We note here that
$\alpha_i\neq 0$ and $\gamma_1\neq 0$ because $\varphi$ is
injective.

We will now compose $\varphi$ with a suitable inner automorphism
to get a nice representative of the class of $\varphi$ in ${\rm
Out}^K(C_r)$ by eliminating $\beta_i$, $i=1,2,3,4$.

Let $y$ be an arbitrary invertible element in $C_r$. Then $y$ is
of the form :
$$y=l_1e_1+l_2e_2+l_3e_3+z,$$
where $l_1,l_2,l_3\in k^*$ and $z\in {\rm rad}\,
C_r$ is a linear combination of the remaining paths of strictly
positive length. Then $y^{-1}$ is easily computed from $yy^{-1}=1.$ Direct
computation gives us that
$$ycy^{-1}=c.$$

Let  $x:=l_1e_1+l_2e_2+l_3e_3+l_4b_1a_1+l_5a_2b_2$, where $l_1$,
$l_2$ and $l_3$ are invertible, and where we set
$l_4:=l_2\alpha_2^{-1}\beta_2$, and
$l_5:=l_2\alpha_3^{-1}\beta_3$. The inner automorphism given by
$x$ has the following action on a set of generators of $C_r$:
$$\begin{aligned}
  xa_1x^{-1}&= l_1l_2^{-1}a_1+l_1l_2^{-2}l_5 a_1a_2b_2, \\
  xa_2x^{-1}&= l_2l_3^{-1} a_2+  l_4l_3^{-1} b_1a_1a_2,\\
  xb_1x^{-1}&= l_2l_1^{-1} b_1+l_1^{-1}l_5a_2b_2b_1,\\
  xb_2x^{-1}&= l_3l_2^{-1} b_2+l_3l_2^{-2}l_4 b_2b_1a_1,\\
\end{aligned}$$
\begin{center}$xcx^{-1}=c,\quad xe_ix^{-1}=e_i,\,\,\ i=1,2,3.$\end{center}

\noindent We denote  by $f^x$ the inner automorphism given by this
specific $x$, and we define $\varphi_1:=f^x\circ \varphi$. This is
an element of ${\rm Out}^K(C_r)$ that is a nice class
representative. Its action on our set of generators is given by
$$\begin{aligned}
  \varphi_1(a_1)&= l_1l_2^{-1}\alpha_1 a_1+(\alpha_1l_1l_2^{-2}l_5+\beta_1l_1l_2^{-1} )a_1a_2b_2, \\
  \varphi_1(a_2)&= l_2l_3^{-1}\alpha_2 a_2+ (\alpha_2 l_4l_3^{-1}+\beta_2l_2l_3^{-1}) b_1a_1a_2,\\
  \varphi_1(b_1)&= l_2l_1^{-1}\alpha_3 b_1+(\alpha_3l_1^{-1}l_5+\beta_3l_2l_1^{-1})a_2b_2b_1,\\
  \varphi_1(b_2)&= l_3l_2^{-1}\alpha_4 b_2+(\alpha_4 l_3l_2^{-2}l_4+\beta_4l_3l_2^{-1}) b_2b_1a_1,\\
\varphi_1(e_i)&=e_i,\,\,\ i=1,2,3,\\
\varphi_1(c)&=c.\end{aligned}$$ We have chosen $l_4$ and $l_5$ in
such a way that, in the above equations, the coefficients of the
paths of length 3 are all equal to 0. The automorphism
$\phi:=f^w\circ \varphi_1$, where $f^w$ is the inner automorphism
given by $w:=l_1^{-1}e_1+l_2^{-1}e_2+l_3^{-1}e_3$, represents the
same class in ${\rm Out}^K(C_r)$ as $\varphi$. It has the
following action on a set of algebra generators:
$$\begin{aligned}
  \phi(e_i)&= e_i, \quad i=1,2,3,\\
  \phi(a_1)&= \alpha_1 a_1, \\
  \phi(a_2)&= \alpha_2a_2,\\
  \phi(b_1)&= \alpha_3b_1, \\
  \phi(b_2)&= \alpha_4b_2, \\
  \phi(c)&= \sum_{i=1}^r\gamma_ic^i.
\end{aligned}$$

\noindent We see that the
$(r+4)$-tuple
$(\alpha_1,\alpha_2,\alpha_3,\alpha_4,\gamma_1,\dots,\gamma_r)$  completely determines $\phi$,
where $\alpha_i,\,i=1,2,3,4,$ and $\gamma_1$ belong to $k^*$ and
$\gamma_2,\dots,\gamma_r\in k$. From the relations of $C_r$ we
have that $\alpha_1\alpha_2\alpha_3\alpha_4=\gamma_1^r.$ It
follows that an arbitrary element  $\phi$ of ${\rm Out}^K(C_r)$ is
determined by an $(r+3)$-tuple, say
$(\alpha_1,\alpha_2,\alpha_3,\gamma_1,\dots,\gamma_r)$, where
$\alpha_4=(\alpha_1\alpha_2\alpha_3)^{-1}\gamma_1^r$. Composition
of homomorphisms induces a group operation on the set of
$(r+3)$-tuples, i.e.\ on the set $k^*\times k^*\times k^*\times
(k^*\times \underbrace{k\times \cdots \times k}_{r-1})$. This is
componentwise multiplication on the first three coordinates and
the operation $*$ of the group $H_r$ from Definition \ref{de:
grupa} on the remaining $r$ coordinates. In other words, we have
the group $(k^*)^3\times H_r$.

Any $(r+3)$-tuple
$(\alpha_1,\alpha_2,\alpha_3,\gamma_1,\dots,\gamma_r)$ gives rise
to a representative of an element of ${\rm Out}^K(C_r)$, i.e.\ we
have an epimorphism from $(k^*)^3\times H_r$ onto ${\rm
Out}^K(C_r)$. The above $(r+3)$-tuple gives us the same class in
${\rm Out}^K(C_r)$ as the $(r+3)$-tuple
$(l_1l_2^{-1}\alpha_1,l_2l_1^{-1}\alpha_2,l_2l_3^{-1}\alpha_3,\gamma_1,\dots,\gamma_r),$
where $l_1,l_2$ and $l_3$ are arbitrary elements from $k^*$. This
corresponds to multiplication by an inner automorphism given by
$l_1e_1+l_2e_2+l_3e_3$. If we set $l_1l_2^{-1}=w$, and 
$l_2l_3^{-1}=v$, then $(k^*)^3\times H_r/R$, where $R$ is the
subgroup generated by all $(r+3)$-tuples of the form
$(w,w^{-1},v,1,0,\dots,0),$ where $v,w\in k^*$, is isomorphic to
${\rm Out}^K(C_r)$. This quotient is isomorphic to the direct
product of one copy of the multiplicative group $k^*$ and a copy
of the group $H_r$. Thus, we see that ${\rm Out}^K(C_r)$ is a
connected algebraic group, and it follows that it is equal to
${\rm Out}^0(C_r)$.

\begin{te}
Let $A$ be one of the algebras  $A_r$, $B_r$ or $C_r$. Then
$${\rm Out}^0(A)\cong k^*\times H_r.$$ The maximal tori in ${\rm Out}^0(A)$
are isomorphic to $\mathbf{G}_m\times
\mathbf{G}_m.$
\end{te}
\noindent{\bf Proof.} This follows from the above discussion and the
fact that ${\rm Out}^0(A)$ is preserved under derived equivalence.
$\blacksquare$

\begin{co}\label{coGradings}
Let $A$ be one of the algebras  $A_r$, $B_r$ or $C_r$. Let $T$ be
a maximal torus in ${\rm Out}(A)$. Then up to graded Morita
equivalence the gradings on $A$ are in one-to-one correspondence
with conjugacy classes in ${\rm Out}(A)$ of cocharacters of ${\rm
Out}(A)$ whose image is in $T$. Up to graded Morita equivalence
the gradings on $A$ are parameterized by the corresponding pairs
of integers.
\end{co}
\noindent{\bf Proof.} By Proposition \ref{cokarak}, up to graded
Morita equivalence the gradings on $A$ are given by conjugacy
classes in ${\rm Out}(A)$ of the algebraic group homomorphisms
from $\mathbf{G}_m$ to ${\rm Out}(A)$.  Let $T^{\prime}$ be
another maximal torus in ${\rm Out}(A)$ and let $f$ be a
cocharacter of ${\rm Out}(A)$ such that its image is contained in
$T^{\prime}$. Since any two maximal tori in ${\rm Out}(A)$ are
conjugate, there exists an invertible element $a$ such that
$aT^{\prime}a^{-1}=T$. The cocharacter given by $x\mapsto
af(x)a^{-1}$, $x\in \mathbf{G}_m$, is conjugate to $f$ and its
image is contained in $T$. This cocharacter gives rise to a
grading which is graded Morita equivalent to the grading given by
$f$. It follows that when classifying gradings on $A$ up to graded
Morita equivalence it is  sufficient to consider cocharacters
whose image is in $T$. Algebraic group homomorphisms from
$\mathbf{G}_m$ to $T\cong\mathbf{G}_m\times \mathbf{G}_m$ are in
one-to-one correspondence with $\mathbb{Z}^2$. $\blacksquare$

\begin{co}
Up to graded Morita equivalence the gradings on $C_r$, $r\geq 2$,
are in  one-to-one correspondence with $\mathbb{Z}^2$.
\end{co}
\noindent {\bf Proof.} From the relations of $C_r$ it follows that
${\rm Out}(C_r)={\rm Out}^K(C_r)$. Let $T$ be the maximal torus in
${\rm Out}(C_r)$ consisting of the $(r+1)$-tuples of the form
$(v,d_1,0,\dots,0)$, where $v,d_1\in k^*$.  Let $\pi_1$ and
$\pi_2$ be  the cocharacters of $T$ corresponding to the pairs of
integers $(m_1,m_2)$ and $(n_1,n_2)$ respectively. If $\pi_1$ and
$\pi_2$ are conjugate in ${\rm Out}(C_r)$, then from the
multiplication in ${\rm Out}(C_r)$ it follows that $m_1=n_1$ and
$m_2=n_2$. $\blacksquare$

\begin{rem}{\rm
There are cases where the group of outer automorphisms of a given
algebra $A$ strictly contains the group of outer automorphisms
fixing the isomorphism classes of simple modules. In this case it
is possible that $N_{{\rm Out}(A)}(T)$ is not contained in ${\rm
Out}^0(A)$, where $T$ is a maximal torus in ${\rm Out}(A)$.

For example, for the remaining two families $A_r$ and $B_r$, the
group of outer automorphisms strictly contains the group of outer
automorphisms fixing the isomorphism classes of simple modules.
This is because there are outer automorphisms in ${\rm Out}(A)$,
where $A$ is $A_r$ or $B_r$, that interchange $e_2$ and $e_3$, and
fix $e_1$. Also, ${\rm Out}^K(A)$ is not necessarily connected,
i.e.\ it is not equal to ${\rm Out}^0(A)$. In this case $N_{{\rm
Out}(A)}(T)$ is not contained in ${\rm Out}^0(A)$, and for
different pairs of integers we get gradings that are graded Morita
equivalent. Thus, $A_r$ and $C_r$ are derived equivalent, but
$N_{{\rm Out}(A_r)}(T)\ncong N_{{\rm Out}(C_r)}(T^{\prime})$,
where $T$ and $T^{\prime}$ are maximal tori.

This tells us that derived equivalent algebras, in general, do not
have the same number of gradings up to graded Morita equivalence.}
\end{rem}

\subsection{Transfer of gradings via derived equivalences}We
will use derived equivalences between $A_r$, $B_r$ and $C_r$ to
transfer gradings from $A_r$ to $B_r$ and $C_r$. The tilting
complexes that we use in this section have been constructed by
Linckelmann in [\ref{Link}].

We assume that $A_r$ is graded in such a way that the vertices and
the arrows of the quiver of $A_r$ are homogeneous. Moreover, we
assume that ${\rm deg}(a_1)=\alpha_1$, ${\rm deg}(a_2)=\alpha_2$,
${\rm deg}(b_1)=\beta_1$, ${\rm deg}(b_2)=\beta_2$ and ${\rm
deg}(c)=\sigma$.  We set
$\Sigma:=\alpha_1+\alpha_2+\beta_2+\beta_2$.

By Example \ref{ex:quiverGrading}, the graded radical layers of
the projective indecomposable  $A_r$-modules  with respect to this
grading are: 

{\small
$$
\begin{array}{ccccc}
&&S_1&&0\\
\alpha_2&S_2&&S_3&\beta_2\\
\alpha_1+\alpha_2&S_1&&S_1&\beta_1+\beta_2\\
\alpha_1+\alpha_2+\beta_2&S_3&&S_2&\beta_1+\beta_2+\alpha_2\\
\Sigma&S_1&&S_1&\Sigma\\
\vdots&&\vdots&&\vdots\\
(r-1)\Sigma+\alpha_2&S_2&&S_3&(r-1)\Sigma+\beta_2\\
r\Sigma-\beta_1-\beta_2&S_1&&S_1&r\Sigma-\alpha_1-\alpha_2\\
r\Sigma-\beta_1&S_3&&S_2&r\Sigma-\alpha_1\\
&&S_1&&r\Sigma
\end{array},$$
$\phantom{a}$
$$\begin{array}{cc}
S_2&0\\
S_1&\alpha_1\\
S_3&\alpha_1+\beta_2\\
S_1&\alpha_1+\beta_1+\beta_2\\
S_2&\Sigma\\
\vdots\\
S_1&(r-1)\Sigma+\alpha_1\\
S_3&r\Sigma-\alpha_2-\beta_1\\
S_1&r\Sigma-\alpha_2\\
S_2&r\Sigma
\end{array},\quad \quad
\begin{array}{cc}
S_3&0\\
S_1&\beta_1\\
S_2&\beta_1+\alpha_2\\
S_1&\beta_1+\alpha_1+\alpha_2\\
S_3&\Sigma\\
\vdots\\
S_1&(r-1)\Sigma+\beta_1\\
S_2&r\Sigma-\alpha_1-\beta_2\\
S_1&r\Sigma-\beta_2\\
S_3&r\Sigma
\end{array}.
$$}Here, numbers to the left or right of the composition
factors denote degrees of the corresponding composition factors.

Let $T_1$ be the complex given by $\, T_1\, :
\xymatrix{P_2\langle-\alpha_2\rangle\oplus
P_3\langle-\beta_2\rangle\ar[rr]^{\,\quad\quad\quad(\gamma_2,
\delta_2)}&& P_1}$, where $P_1$ is in degree 1, and $\gamma_2,
\delta_2$ are given by right multiplication by $a_2$ and $b_2$
respectively. Let $T_2$  and $T_3$ be the stalk complexes with
$P_2$ and $P_3$ respectively in degree 0. A complex $T$ that tilts
from $A_r$ to $B_r$ is given by  the direct sum $T:=T_1\oplus T_2
\oplus T_3$.

Viewing $T$ as a graded object and calculating ${\rm
Homgr}_{K^b(P_{A_r})}(T,T)$  as a graded vector space will give us
a grading on $B_r$. It is clear that
$$
\begin{aligned}
{\rm Homgr}_{K^b(P_{A_r})}(T_2, T_2)&\cong {\rm
Homgr}_{A_r}(P_2,P_2)\cong \bigoplus_{t=0}^r k\langle
-t\Sigma\rangle,\\
{\rm Homgr}_{K^b(P_{A_r})}(T_3, T_3)&\cong  {\rm
Homgr}_{A_r}(P_3,P_3)\cong \bigoplus_{t=0}^r k\langle
-t\Sigma\rangle,\\
{\rm Homgr}_{K^b(P_{A_r})}(T_2, T_3)&\cong {\rm
Homgr}_{A_r}(P_2,P_3)  \cong \bigoplus_{t=0}^{r-1} k\langle
-(\beta_1+\alpha_2)-t\Sigma \rangle,\\
{\rm Homgr}_{K^b(P_{A_r})}(T_3, T_2)&\cong {\rm
Homgr}_{A_r}(P_3,P_2)  \cong \bigoplus_{t=0}^{r-1} k\langle
-(\beta_2+\alpha_1)-t\Sigma \rangle.\\
\end{aligned}
$$

\noindent It follows that  ${\rm deg}(d_2)=\alpha_1+\beta_2$ and
${\rm deg}(c_2)=\beta_1+\alpha_2$ in the quiver of $B_r$.

\noindent Also, non-zero maps in   ${\rm
Homgr}_{K^b(P_{A_r})}(T_1, T_2)$ and ${\rm
Homgr}_{K^b(P_{A_r})}(T_1, T_3)$ have to map surjectively
$P_2\oplus P_3$ onto $P_2$ and $P_3$ respectively. We conclude that  ${\rm
Homgr}_{K^b(P_{A_r})}(T_1, T_2)\cong k\langle \alpha_2 \rangle$,
and ${\rm Homgr}_{K^b(P_{A_r})}(T_1, T_3)\cong k\langle \beta_2
\rangle$. It follows that ${\rm deg}(c_1)=-\alpha_2$ and ${\rm
deg}(d_3)=-\beta$ in the quiver of $B_r$.

Every non-zero map in ${\rm Homgr}_{K^b(P_{A_r})}(T_2, T_1)$ has
to map ${\rm top}\, P_2$ onto  ${\rm soc}\, P_2$. It follows that
${\rm Homgr}_{K^b(P_{A_r})}(T_2, T_1)\cong k\langle
-\alpha_2-r\Sigma\rangle$, and similarly we deduce that ${\rm
Homgr}_{K^b(P_{A_r})}(T_3, T_1)\cong k\langle -\beta_2-r\Sigma
\rangle$.  This implies that ${\rm deg}(c_3)=\beta_2+r\Sigma$ and
${\rm deg}(d_1)=\alpha_2+r\Sigma$.

From the above computation  we get a grading on $B_r$. With respect
to this grading, the graded quiver of $B_r$ is given by
$$
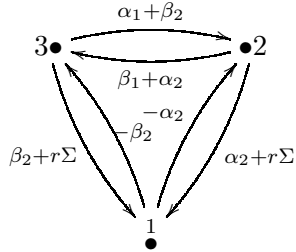
\xymatrix{3\bullet\ar@{<-}@/^/[ddr]^(.55){\!\!\!\!-\beta_2}\ar@/^/[rr]^{\alpha_1+\beta_2}&&\bullet 2\ar@/^/[ll]^{\beta_1+\alpha_2}\ar@/^/[ddl]^{\alpha_2+r\Sigma}\\
&&\\
&\stackrel{1}{\bullet}\ar@/^/[uur]^(.55){-\alpha_2}\ar@{<-}@/^/[uul]^{\beta_2+r\Sigma}&}
$$

If we assume that we started with the tight grading on $A_r$, i.e.
if we assume that the arrows of the quiver of $A_r$ are in degree
1, then the resulting graded quiver of $B_r$ is given by
$$
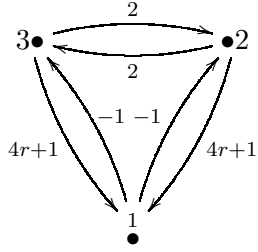
\xymatrix{3\bullet\ar@{<-}@/^/[ddr]^{\!\!-1}\ar@/^/[rr]^{2}&&\bullet 2\ar@/^/[ll]^{2}\ar@/^/[ddl]^{4r+1}\\
&&\\
&\stackrel{1}{\bullet}\ar@/^/[uur]^{\,-1}\ar@{<-}@/^/[uul]^{4r+1}&}
$$

We remark here that the resulting grading on $B_r$ is not tight.
Moreover, it is not a positive grading.  This example tells us
that tightness and positivity of a grading are not preserved under
derived equivalences.  We state this
 known fact in the following proposition.
\begin{prop}
Tightness and positivity of a grading are not preserved, in
general, under the transfer of gradings via derived equivalence.
\end{prop}

Let us now assume that the algebra $B_r$ is graded in such a way
that the vertices and the arrows of the quiver of $B_r$ are
homogeneous. Furthermore, we assume that ${\rm
deg}(c_1)=\gamma_1$, ${\rm deg}(c_2)=\gamma_2$, ${\rm
deg}(c_3)=\gamma_3$, ${\rm deg}(d_1)=\delta_1$, ${\rm
deg}(d_2)=\delta_2$ and ${\rm deg}(d_3)=\delta_3$.  We set
$\Sigma:=\gamma_2+\delta_2$.

The graded radical layers of the projective indecomposable
$B_r$-modules are: 

{\small
$$
\begin{array}{ccccc}
&&S_1&&0\\
\delta_1&S_2&&S_3&\gamma_3\\
&&S_1&&r\Sigma
\end{array},$$ 

$\phantom{a}$

$$
\begin{array}{ccccc}
&&S_2&&0\\
&&&S_3&\delta_2\\
&&&S_2&\Sigma\\
\gamma_1&S_1&&S_3&\Sigma+\delta_2\\
&&&S_2&2\Sigma\\
&&&\vdots&\\
&&&S_3&r\Sigma-\gamma_2\\
&&S_2&&r\Sigma
\end{array},\quad
\begin{array}{ccccc}
&&S_3&&0\\
&&&S_2&\gamma_2\\
&&&S_3&\Sigma\\
\delta_3&S_1&&S_2&\Sigma+\gamma_2\\
&&&S_3&2\Sigma\\
&&&\vdots&\\
&&&S_2&r\Sigma-\delta_2\\
&&S_3&&r\Sigma
\end{array}.
$$}


We will now transfer this grading from $B_r$ to $C_r$. Let $T_1$
and $T_3$  be the stalk complexes with $P_1$ and $P_3$
respectively in degree 0. Let $T_2$ be the complex $$T_2\, :\,\,
\xymatrix{P_1\langle -\gamma_1\rangle\oplus P_3\langle -\delta_2
\rangle \ar[rr]^{\hspace{15mm}(\rho_1 ,\tau_2)}&& P_2},$$ where
$P_2$ is in degree 1, and $\rho_1, \tau_2$ are given by right
multiplication by $c_1$ and $d_2$  respectively. Define $T$ to be
the direct sum  $T:=T_1\oplus T_2\oplus T_3$. The complex $T$  is
a tilting complex for $B_r$  and ${\rm End}_{K^b(P_{B_r})}(T)\cong
C_r^{op}$.

 As above, we conclude that the space ${\rm
Homgr}_{K^b(P_{B_r})}(T_3,T_3)$ is isomorphic to
$\bigoplus_{t=0}^r k\langle -t\Sigma\rangle$, ${\rm
Homgr}_{K^b(P_{B_r})}(T_1,T_1)$ is isomorphic to  $k\langle 0
\rangle\oplus k\langle -r\Sigma\rangle$,  ${\rm
Homgr}_{K^b(P_{B_r})}(T_3,T_1)\cong k\langle -\gamma_3\rangle$,
and ${\rm Homgr}_{K^b(P_{B_r})}(T_1,T_3)\cong k\langle -\delta_3
\rangle$. It follows that ${\rm deg}(c)=\Sigma$ in the quiver of
$C_r$. Since ${\rm ker}(\rho_1,\tau_2)$ contains two copies of
$S_1$, one copy in degree $\delta_3+\delta_2$ and  one copy in
degree $\gamma_1+r\Sigma$, ${\rm
Homgr}_{K^b(P_{B_r})}(T_1,T_2)\cong k\langle -(\delta_3+\delta_2)
\rangle\oplus k\langle -(\gamma_1+r\Sigma) \rangle$. The same arguments give us that  ${\rm
Homgr}_{K^b(P_{B_r})}(T_3,T_2)\cong k\langle -(\gamma_1+\gamma_3)
\rangle\oplus k\langle -(\delta_2+r\Sigma)\rangle.$ Similarly, there are isomorphisms ${\rm Homgr}_{K^b(P_{B_r})}(T_2, T_3)\cong k\langle \delta_2
\rangle\oplus k\langle \gamma_1-\delta_3 \rangle$, and ${\rm
Homgr}_{K^b(P_{B_r})}(T_2, T_1)\cong k\langle \gamma_1
\rangle\oplus k\langle \delta_2-\gamma_3 \rangle$.

Using these data and looking at the relations of  $C_r$, we have
that in the quiver of $C_r$, ${\rm deg}(a_1)=\delta_2+\delta_3$,
${\rm deg}(a_2)=-\delta_2$, ${\rm deg}(b_1)=-\gamma_1$ and ${\rm
deg}(b_2)=\gamma_1+\gamma_3$. With respect to this grading, the
graded quiver of $C_r$ is given by
 $$
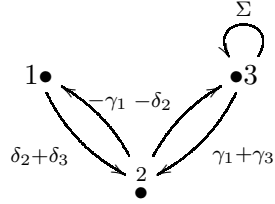
\xymatrix{1\bullet\ar@/_/[dr]_{\delta_2+\delta_3}&&\bullet 3\ar@/^/[dl]^{\gamma_1+\gamma_3}\ar@(ur,ul)_{\Sigma}\\
&\stackrel{2}{\bullet}\ar@/^/[ur]^{\,\,\,\,\,-\delta_2}\ar@/_/[ul]_{\!\!\!\!\!\!-\gamma_1}&}
$$


\noindent The graded radical layers of the projective
indecomposable $C_r$-modules are:

 {\small
$$
\begin{array}{cc}
S_1&\\
S_2&-\gamma_1\\
S_3&\gamma_3\\
S_2&\gamma_3-\delta_2\\
S_1&r\Sigma
\end{array},\quad
\begin{array}{ccccc}
&&S_2&&\\
\delta_2+\delta_3&S_1&&S_3&\gamma_1+\gamma_3\\
\delta_2+\delta_3-\gamma_1&S_2&&S_2&\gamma_1+\gamma_3-\delta_2\\
\delta_2+\delta_3+\gamma_3&S_3&&S_1&\gamma_1+\gamma_3+\delta_3\\
&&S_2&&r\Sigma
\end{array},$$ $$
\begin{array}{ccccc}
&&S_3&&\\
\Sigma&S_3&&S_2&-\delta_2\\
2\Sigma&S_3&&S_1&\delta_3\\
&\vdots&&S_2&\delta_3-\gamma_1\\
(r-1)\Sigma&S_3&&&\\
&&S_3&&r\Sigma
\end{array}.
$$}



\subsection{Positivity and tightness}
\begin{prop}
Let $A$ be one of the algebras $A_r$, $B_r$ and $C_r$. Then $A$
can be positively graded.
\end{prop}
\noindent {\bf Proof.} This follows directly from the relations of
these algebras. For $A_r$ we can set that every arrow is in degree
1 and we will get homogeneous relations. For the algebra $B_r$, if  ${\rm deg}(c_1)={\rm
deg}(d_1)={\rm deg}(c_3)={\rm deg}(d_3)=r$ and ${\rm
deg}(c_2)={\rm deg}(d_2)=1$, then the relations of
$B_r$ are homogeneous. If
${\rm deg}(c)=1$, ${\rm deg}(a_1)={\rm deg}(a_2)={\rm
deg}(b_1)={\rm deg}(b_2)=r$,  then the relations of $C_r$ are homogeneous.  $\blacksquare$

\begin{prop}
For every positive integer $r$, $A_r$ is a tightly graded algebra.
\end{prop}
\noindent {\bf Proof.} From the proof of the previous
proposition, if the vertices of the quiver of $A_r$ are in degree 0, and the arrows are in degree 1, then the ideal of relations of $A_r$ is homogeneous. Therefore, there exists a positive grading on $A$ such that the subalgebra of degree 0 elements is semisimple, and $A$ is generated by the homogeneous elements of degrees 0 and 1. $\blacksquare$

\begin{prop}
The algebra $B_r$ is tightly graded if and only if $r=1$.
\end{prop}
\noindent {\bf Proof.} It is clear that $B_1$ is tightly graded.
Let us assume that $B_r$ is tightly graded. By Lemma \ref{lm:
tight}, for each arrow $a$ of the quiver of $B_r$, there exists a
degree 1 element of the form $a+\sum_i\lambda_iz_i$, where $z_i\in
{\rm rad}^2\,A$  is a path with the same source and the same
target as $a$. It follows that $c_1$, $c_3$, $d_1$ and $d_3$ are
homogeneous elements of degree 1, since there are no other paths
with the same source and the same target. Also, there are degree 1
elements of the form
$$t_{c_2}:=c_2+\sum_{i=1}^{r-1}\lambda_ic_2(d_2c_2)^i,$$
$$t_{d_2}:=d_2+\sum_{i=1}^{r-1}\mu_id_2(c_2d_2)^i,$$
where the $\lambda$'s and $\mu$'s are scalars.

\noindent It follows that $(t_{c_2}t_{d_2})^r=(c_2d_2)^r$ is a
homogeneous element of degree $2r$. Since $(c_2d_2)^r=d_1c_1$ is a
homogeneous element of degree 2, it follows that $r=1$.
$\blacksquare$

\begin{prop}
The algebra $C_r$ is tightly graded if and only if $r=1$ or $r=4$.
\end{prop}
\noindent {\bf Proof.} If $r=1$ or $r=4$, then it is obvious that
$C_r$ is tightly graded.

Let us assume that $C_r$, $r\geq 2$, is tightly graded. By Lemma
\ref{lm: tight}, there are degree 1 elements of the form
$$
\begin{array}{c@{:=}l}
t_{a_1}\,&a_1+\lambda_1a_1a_2b_2,\\
t_{a_2}\,&a_2+\lambda_2b_1a_1a_2,\\
t_{b_1}\,&b_1+\lambda_3a_2b_2b_1,\\
t_{b_2}\,&b_2+\lambda_4b_2b_1a_1,\\
t_{c}\,&c+\sum_{i=2}^{r}\mu_ic^i,
\end{array}
$$
where the $\lambda$'s and $\mu$'s are scalars.

It follows that $b_2b_1a_1a_2=t_{b_2}t_{b_1}t_{a_1}t_{a_2}$ is a
homogeneous element of degree 4. At the same time
$b_2b_1a_1a_2=c^r=t_c^r$ is a homogeneous element of degree $r$.
It follows that $r=4.$ $\blacksquare$

We note here that from the previous propositions it follows that
the existence of a tight grading is not preserved under derived
equivalence, unlike under Morita equivalence (see Proposition
4.4 in [\ref{CPS}]).

It is worth noting that for dihedral blocks with three simple
modules, in every derived equivalence class there is at least one
block that is positively graded and there is at least one block
that is tightly graded. The same statement does not hold for all
derived equivalence classes of tame blocks.

\section{Two simple modules}
Any block with a dihedral defect group and two isomorphism classes
of simple modules is Morita equivalent to some algebra from the
following list (cf.\ [\ref{ErK}] or [\ref{H96}]).
\begin{itemize}
 \item[(1)] For any $r\geq 1$ and $c\in\{0,1\}$ let $D(2A)^{r,c}$ be the algebra defined by the quiver and relations
\begin{multicols}{2}
$
\xymatrix{&0\bullet\ar@(dl,ul)^{\alpha}\ar@/^/[rr]^{\beta}&&\bullet
1\ar@/^/[ll]^{\gamma}}
$

$\!\!\!\!\!\!\gamma\beta=0,$ $\alpha^2=c(\alpha\beta\gamma)^r$,\\
$(\alpha\beta\gamma)^r=(\beta\gamma\alpha)^r$.
\end{multicols}
\item[(2)]For any $r\geq 1$ and $c\in\{0,1\}$ let $D(2B)^{r,c}$
be the algebra defined by the quiver and relations
\begin{multicols}{2}
$
\xymatrix{&0\bullet\ar@(dl,ul)^{\alpha}\ar@/^/[rr]^{\beta}&&\bullet
1\ar@/^/[ll]^{\gamma}\ar@(dr,ur)_{\eta}}
\phantom{blaaaaa}
$

$\!\!\!\!\!\!\!\beta\eta=\eta\gamma=\gamma\beta=0$,\\
$\alpha\beta\gamma=\beta\gamma\alpha$,\\
$\alpha^2=c(\alpha\beta\gamma)$,
$\gamma\alpha\beta=\eta^r$.
\end{multicols}
\end{itemize}

\subsection{Classification of gradings}In [\ref{H96}], Holm proved that for fixed $r$ and $c$, the
algebras $D(2A)^{r,c}$ and $D(2B)^{r,c}$ are derived equivalent.
Since the identity component of the group of outer automorphisms
is invariant under derived equivalence, it is sufficient to
compute this group for $D(2B)^{r,c}$.

As before, for an arbitrary outer automorphism $\varphi$ in
${\rm Out}^K(D(2B)^{r,c})$, we will find a suitable automorphism
that represents the same element as $\varphi$, but which is easy
to work with.

We assume that $\varphi(e_i)=e_i$, for $i=1,2.$ It follows that
$$\begin{aligned}
\varphi(\alpha)&= a_1\alpha+a_2\beta\gamma+a_3\alpha\beta\gamma,\\
\varphi(\beta)&= b_1\beta+b_2\alpha\beta,\\
\varphi(\gamma)&= c_1\gamma+c_2\gamma\alpha,\\
\varphi(\eta)&= \sum_{i=1}^rd_i\eta^i,
  \end{aligned}
$$
for some $a_i,b_i,c_i,d_i\in k$. From the relation $\gamma\beta=0$
we get that $b_1c_2=b_2c_1$.  From the relation
$\eta^r=\gamma\alpha\beta$ it follows that $d_1^r=a_1b_1c_1$.
Since $\varphi(\eta^r)\neq 0$, it follows that $d_1\neq 0$. Hence,
$a_1, b_1$ and $c_1$ are all non-zero. The inner automorphism
given by $y$, where $y:=l_1e_1+l_2e_2+l_3\alpha$ and
$l_3:=l_1c_1^{-1}c_2$, when composed with $\varphi$ has the
following action on a set of generators:
$$\begin{aligned}
y\varphi(e_i)y^{-1}&= e_i,\\
   y\varphi(\eta)y^{-1}&= \sum_{i=1}^rd_i\eta^i,\\
y\varphi(\gamma)y^{-1}&= c_1l_2l_1^{-1}\gamma,\\
y\varphi(\beta)y^{-1}&= b_1l_1l_2^{-1}\beta,\\
y\varphi(\alpha)y^{-1}&=
a_1\alpha+a_2\beta\gamma+a_3\alpha\beta\gamma.
  \end{aligned}
$$
Let $\phi$ be the composition of $y\varphi y^{-1}$ and the inner
automorphism given by $l_1^{-1}e_1+l_2^{-1}e_2$. Then $\phi$
represents the same element in ${\rm Out}^K(D(2B)^{r,c})$ as
$\varphi$.  Its action is given by
$$\begin{aligned}
\phi(e_i)&= e_i,\\
   \phi(\eta)&= \sum_{i=1}^rd_i\eta^i,\\
\phi(\gamma)&= c_1\gamma,\\
\phi(\beta)&= b_1\beta,\\
\phi(\alpha)&= a_1\alpha+a_2\beta\gamma+a_3\alpha\beta\gamma.
  \end{aligned}
$$
It follows that an arbitrary automorphism in ${\rm
Out}^K(D(2B)^{r,c})$ is completely determined by an $(r+5)$-tuple
$(a_1,a_2,a_3,b_1,c_1, d_1,\dots, d_r)$. By an elementary, but a
tedious calculation, one can show that it is not possible to
eliminate coefficients $a_2$ and $a_3$ by composing $\phi$ with
inner automorphisms.

 We have a map from the set of all $(r+5)$-tuples onto ${\rm
Out}^K(D(2B)^{r,c})$. Composition  of morphisms gives us the group
multiplication on the set of all $(r+5)$-tuples.

From $d_1^r=a_1b_1c_1$ it follows that one of these four
coefficients, say $a_1$, is determined by the remaining three.

If $c=0$, then there are no further restrictions to the
coefficients of $\varphi$. In this case,  $\varphi$ is determined
by the $(r+4)$-tuple $(a_2,a_3,b_1,c_1,d_1,\dots,d_r)$, where
$b_1,c_1,d_1\in k^*$. The multiplication of these $(r+4)$-tuples
is given by composition of the corresponding automorphisms, where
we replace $a_1$ with $d_1^r(b_1c_1)^{-1}$. If $G$ is the group of
all such $(r+4)$-tuples, then the multiplication is given by:
$$(a_2^{\prime},a_3^{\prime},b_1^{\prime},c_1^{\prime},\mathbf{d^{\prime}})*(a_2,a_3,b_1,c_1,\mathbf{d})=$$
$$=(d_1^r(b_1c_1)^{-1}a_2^{\prime}+a_2b_1^{\prime}c_1^{\prime},d_1^r(b_1c_1)^{-1}a_3^{\prime}+a_3(d_1^{\prime})^r,b_1b_1^{\prime},c_1c_1^{\prime},\mathbf{d}\mathbf{d^{\prime}}),$$
where $\mathbf{d}=(d_1,\dots,d_r)$ and
$\mathbf{d^{\prime}}=(d_1^{\prime},\dots,d_r^{\prime})$, and the
product $\mathbf{d}\mathbf{d^{\prime}}$ is the product of elements
of the group $H_r$ from Definition \ref{de: grupa}.

Thus, we have a map from  the group $G$ of all $(r+4)$-tuples
onto  the group ${\rm Out}^K(D(2B)^{r,c})$. The kernel of this epimorphism
is given by the $(r+4)$-tuples that correspond to inner
automorphisms. Let $R$ be the subgroup of $G$ generated by all
$(r+4)$-tuples that correspond to inner automorphisms. The
$(r+4)$-tuple $(a_2,a_3,b_1,c_1,\mathbf{d})$ represents the same
class in the quotient group $M:=G/R$ as
$(a_2,a_3,l_1l_2^{-1}b_1,l_1^{-1}l_2c_1,\mathbf{d})$, where
$l_1,l_2\in k^*$. In particular, if $l_1l_2^{-1}=c_1$, then the
$(r+4)$-tuple $(a_2,a_3,b_1,c_1,\mathbf{d})$ represents the same
element as the $(r+4)$-tuple $(a_2,a_3,b_1c_1,1,\mathbf{d})$. If
$v=b_1c_1$, then $M$ can be seen as the group consisting of $(r+3)$-tuples $(a_2,a_3,v,\mathbf{d})$, where the multiplication
is defined by:
$$(a_2^{\prime},a_3^{\prime},v^{\prime},\mathbf{d^{\prime}})*(a_2,a_3,v,\mathbf{d})=(d_1^rv^{-1}a_2^{\prime}+a_2v^{\prime},d_1^rv^{-1}a_3^{\prime}+a_3(d_1^{\prime})^r,vv^{\prime},\mathbf{d}\mathbf{d^{\prime}}).$$

\begin{prop}
Let $M$ be as above and let $A$ be $D(2B)^{r,0}$ or $D(2A)^{r,0}$.
There is an isomorphism of  groups $${\rm Out}^0(A)\cong M.$$ The
maximal tori in ${\rm Out}^0(A)$ are isomorphic to
$\mathbf{G}_m\times \mathbf{G}_m.$
\end{prop}
\noindent {\bf Proof.} From the above discussion follows that
${\rm Out}^K(D(2B)^{r,c})$ is isomorphic to $M$. Because ${\rm
Out}^K(D(2B)^{r,c})$ is connected, it is equal to the identity component ${\rm
Out}^0(D(2B)^{r,c})$. The identity component of the group of outer
automorphisms is invariant under derived equivalence. Hence, the
first statement of the proposition is true.

The subgroup $L$ of $M$ which is generated by the $(r+3)$-tuples
of the form $(a_2,a_3,1,1,d_2,\dots,d_r)$ is a normal subgroup of
$M$. The subgroup $T$ of $M$ generated by the $(r+3)$-tuples of
the form $(0,0,v,d_1,0,\dots,0)$ is isomorphic to the quotient
$M/L$. It follows that $M$ is isomorphic to the semidirect product
$L\rtimes T$. The group $L$ is unipotent and the group $T$ is
semisimple. Since $T\cong \mathbf{G}_m\times \mathbf{G}_m$, it
follows that the maximal tori in ${\rm Out}^K(D(2B)^{r,0})$ are
isomorphic to $\mathbf{G}_m\times \mathbf{G}_m$. $\blacksquare$

\begin{co}
Let $A$ be one of the algebras  $D(2B)^{r,0}$ or $D(2A)^{r,0}$.
Let $T$ be a maximal torus in ${\rm Out}(A)$. Then up to graded
Morita equivalence the gradings on $A$ are in one-to-one
correspondence with conjugacy classes in ${\rm Out}(A)$ of
cocharacters of ${\rm Out}(A)$ whose image is in $T$. Up to graded
Morita equivalence the gradings on $A$ are parameterized by the
corresponding pairs of integers.
\end{co}
\noindent{\bf Proof.} The proof is the same as the proof of
Corollary \ref{coGradings} $\blacksquare$

\begin{co}
Up to graded Morita equivalence the gradings on $D(2B)^{r,0}$ are
in  one-to-one correspondence with $\mathbb{Z}^2$.
\end{co}
\noindent {\bf Proof.} It follows from the relations of $D(2B)^{r,0}$ that an arbitrary outer automorphism has to fix the vertices of
the quiver of $D(2B)^{r,0}$. Hence,  ${\rm Out}(D(2B)^{r,0})={\rm
Out}^K(D(2B)^{r,0})$. Let $T$ be the maximal torus 
consisting of the $(r+4)$-tuples of the form
$(0,0,v,d_1,0,\dots,0)$, where $v,d_1\in k^*$. Let $\pi_1$ and
$\pi_2$ be  the cocharacters of $T$ corresponding to the pairs of
integers $(m_1,m_2)$ and $(n_1,n_2)$ respectively. If $\pi_1$ and
$\pi_2$ are conjugate in ${\rm Out}(D(2B)^{r,0})$, then from the
multiplication in ${\rm Out}(D(2B)^{r,0})$ it follows that
$m_1=n_1$ and $m_2=n_2$. $\blacksquare$

As in the case of three simple modules, the same remarks about the
gradings on $D(2A)^{r,0}$ hold, since ${\rm Out}^K(D(2A)^{r,0})$
is not a connected group.

If $c=1$ there is an additional restriction to the coefficients of
$\varphi$ coming from the relation $\alpha^2=\alpha\beta\gamma$.
From this relation we have that $a_1=b_1c_1$. This implies that
$b_1c_1=\sqrt{d_1^r}$. It follows that one of these coefficients,
say $b_1$, is determined by the remaining two. In this case
$\varphi$ is determined by the $(r+3)$-tuple
$(a_2,a_3,c_1,d_1,\dots,d_r)$. We have a map from the group $G$ of
all $(r+3)$-tuples onto ${\rm Out}^K(D(2B)^{r,1})$. The
multiplication in $G$ is the same as before, in this case we just
set $b_1c_1=\sqrt{d_1^r}$. The kernel of the above map is the
subgroup $R$ generated by all $(r+3)$-tuples corresponding to
inner automorphisms.  It follows that in the quotient group $G/R$,
the $(r+3)$-tuple $(a_2,a_3,c_1,d_1,\dots,d_r)$ represents the
same element as the $(r+3)$-tuple $(a_2,a_3,1,d_1,\dots,d_r)$. We
have that $G/R$ is the group consisting of $(r+2)$-tuples
$(a_2,a_3,d_1,\dots,d_r)$, with the multiplication given by:
$$(a_2^{\prime},a_3^{\prime},\mathbf{d^{\prime}})*(a_2,a_3,\mathbf{d})=(\sqrt{d_1^r}a_2^{\prime}+a_2\sqrt{(d_1^{\prime})^r},\sqrt{d_1^r}a_3^{\prime}+a_3(d_1^{\prime})^r,\mathbf{d^{\prime}}\mathbf{d}).$$

\begin{prop}\label{dvamodulajunik}
Let $A$ be one of the algebras  $D(2B)^{r,1}$ or $D(2A)^{r,1}$.
Let $G$ and $R$ be as above. Then ${\rm Out}^0(A)\cong G/R$. The
maximal tori in ${\rm Out}^0(A)$ are isomorphic to $\mathbf{G}_m$.
Up to graded Morita equivalence and rescaling there is a unique
grading on the algebra $A$.
\end{prop}
\noindent{\bf Proof.} It is obvious that ${\rm
Out}^K(D(2B)^{r,1})$ is connected, hence it is equal to its
identity component ${\rm Out}^0(D(2B)^{r,1})$. That ${\rm
Out}^0(A)\cong G/R$ follows from the above discussion and the fact
that the identity component of the group of outer automorphisms is
invariant under derived equivalence. It is easily verified that
$G/R\cong L\rtimes T$, where $T$ is the subgroup generated by all
$(r+2)$-tuples of the form $(0,0,d_1,0,\dots,0)$, and $L$ is the
subgroup generated by all $(r+2)$-tuples of the form
$(a_2,a_3,1,d_2,\dots,d_r)$. It follows that the maximal tori are
isomorphic to $\mathbf{G}_m.$ By Lemma \ref{KlGr}, there is a unique
grading on $A$ up to graded Morita equivalence and rescaling.
$\blacksquare$

An easy corollary of our results is that for different values of
the scalar $c$ we get algebras that are not derived equivalent.
This statement follows from the fact that ${\rm Out}^0(A)$ is
invariant under derived equivalence. On the other hand, ${\rm
Out}^0(D(2B)^{r,0})$ and ${\rm Out}^0(D(2B)^{r,1})$ are not
isomorphic because they do not have isomorphic maximal tori. Even
though this is known (cf.\ [\ref{H99}], Proposition 3.1), we
record it in the following corollary.

\begin{co}
Let $C^{r,0}$ be one of the algebras $D(2A)^{r,0}$ or
$D(2B)^{r,0}$, and let $C^{r,1}$ be one of the algebras
$D(2A)^{r,1}$ or $D(2B)^{r,1}$. Then $C^{r,0}$ and  $C^{r,1}$ are
not derived equivalent.
\end{co}

\subsection{Transfer of gradings via derived equivalences}  We will
use tilting complexes given in [\ref{H96}] to transfer gradings
from $D(2A)^{r,c}$ to $D(2B)^{r,c}$. Let us fix an integer $r$ and
$c\in \{0,1\}$,  and assume that $D(2A)^{r,c}$ is graded  in such
a way that the vertices and the arrows of the quiver of
$D(2A)^{r,c}$ are homogeneous.  We assume that  the arrows
$\alpha,\beta$ and $\gamma$ of the quiver of $D(2A)^{r,c}$ are in
degrees $d_1,d_2$ and $d_3$ respectively. We set $d:=d_1+d_2+d_3.$

The graded radical layers of the projective indecomposable
$D(2A)^{r,c}$-modules are: 

 {\small $$
\begin{array}{ccccc}
&&S_0&&0\\
d_1&S_0&&S_1&d_3\\
d_1+d_3&S_1&&S_0&d_2+d_3\\
d&S_0&&S_0&d\\
\vdots&\vdots&&\vdots&\vdots\\
(r-1)d+d_1&S_0&&S_1&(r-1)d+d_3\\
rd-d_2&S_1&&S_0&rd-d_1\\
&&S_0&&rd
\end{array},
\quad
\begin{array}{cc}
S_1&0\\
S_0&d_2\\
S_0&d_1+d_2\\
S_1&d\\
\vdots\\
S_0&(r-1)d+d_2\\
S_0&(r-1)d-d_3\\
S_1&rd
\end{array}.
$$}
Since the relations are homogeneous we have that
$(r-2)d_1+rd_2+rd_3=0$ if $c=1$. In this case $d_1,d_2$ and $d_3$
cannot all be non-negative (unless they are all equal to zero). If
$c=0$, all relations are trivially homogeneous and we can choose
$d_1,d_2$ and $d_3$ arbitrarily. In particular, if $c=0$, then
$D(2A)^{r,c}$ is a tightly graded algebra.

A graded tilting complex $T:=T_0\oplus T_1$ of projective
$D(2A)^{r,c}$-modules that tilts from $D(2A)^{r,c}$ to
$D(2B)^{r,c}$  is given by the direct sum of the complex $T_1$,
which is the stalk complex with $P_1$ in degree 0, and the complex
$$T_0\, :\,\, \xymatrix{0\ar[r]&P_1\langle -d_3\rangle\oplus
P_1\langle -(d_1+d_3)\rangle
\ar[rr]^{\quad\quad\,\,\,\,\,\,\hspace{10mm}\,\,(\gamma,
\gamma\alpha)}&& P_0},$$ where $P_0$ is in degree 1, and where
$\gamma$ and $\gamma\alpha$ are given by right multiplication by
$\gamma$ and $\gamma\alpha$ respectively. It was shown in
[\ref{H96}] that $T$ is a tilting complex for $D(2A)^{r,c}$  and
that ${\rm End}_{K^b(P_{D(2A)^{r,c}})}(T)\cong
(D(2B)^{r,c})^{op}$. Viewing $T$ as a graded object and
calculating ${\rm Endgr}_{K^b(P_{D(2A)^{r,c}})}(T)$  as a graded
vector space will give us a grading on $D(2B)^{r,c}$.

From ${\rm Homgr}_{K^b(P_{D(2A)^{r,c}})}(T_1, T_1)\cong
\bigoplus_{t=0}^r k\langle -td\rangle$ we have ${\rm deg
}(\eta)=d$.

To calculate ${\rm Homgr}_{K^b(P_{D(2A)^{r,c}})}(T_1, T_0)$ notice
that this space is isomorphic to ${\rm Homgr}_{D(2A)^{r,c}}(P_1,
{\rm ker}(\gamma,\gamma\alpha))$. Non-zero maps in the latter
space have to map ${\rm top}\, P_1$ to ${\rm soc}\,
P_1\langle-d_3\rangle$, or to ${\rm soc}\, P_1\langle
-(d_1+d_3)\rangle$. This gives us that $${\rm
Homgr}_{K^b(P_{D(2A)^{r,c}})}(T_1, T_0)\cong k\langle
-(rd+d_3)\rangle \oplus k\langle -(rd+d_1+d_3)\rangle.$$ Since the
only non-zero paths in the quiver of $D(2B)^{r,c}$  that start at
vertex $1$ and end at vertex $0$ are $\gamma$ and $\gamma\alpha$,
then $$\{{\rm deg}(\gamma),{\rm deg
}(\gamma\alpha)\}=\{rd+d_3, rd+d_1+d_3\}.$$

To calculate ${\rm Homgr}_{K^b(P_{D(2A)^{r,c}})}(T_0, T_1)$ notice
that non-zero maps in this space have to map
$P_1\langle-d_3\rangle$ or $P_1\langle -(d_1+d_3)\rangle $ onto
$P_1$. It follows that $${\rm Homgr}_{K^b(P_{D(2A)^{r,c}})}(T_0,
T_1)\cong k\langle d_3 \rangle \oplus k\langle d_1+d_3 \rangle.$$
Since the only non-zero paths in the quiver of $D(2B)^{r,c}$ that
start at vertex $0$ and end at vertex $1$ are $\beta$ and
$\alpha\beta$, we have that $$\{{\rm deg}(\beta),{\rm deg
}(\alpha\beta)\}=\{-d_3, -d_1-d_3\}.$$

There are two choices for ${\rm deg}(\alpha)$. If ${\rm
deg}(\alpha)=d_1$, then ${\rm deg}(\beta)=-d_1-d_3$ and ${\rm
deg}(\gamma)=rd+d_3$. This gives us a grading on $D(2B)^{r,c}$. If
${\rm deg}(\alpha)=-d_1$, then ${\rm deg}(\beta)=-d_3$ and ${\rm
deg}(\gamma)=rd+d_1+d_3$. This will not give us a grading on
$D(2B)^{r,c}$ if $c = 1$, because the relations are not
homogeneous. If $c = 0$, this grading is the same as the previous
one via suitable substitution of the integers $d_1, d_2, d_3$,
i.e.\ we get this grading from the former grading if we choose
$-d_1$, $d_1+d_2$, $d_1+d_3$ instead of $d_1$, $d_2$ and $d_3$
respectively for the degrees of the corresponding arrows.

With respect to this resulting grading, the graded
quiver of $D(2B)^{r,c}$  is given by

$$
\xymatrix{0\bullet\ar@(dl,ul)^{d_1}\ar@/^/[rr]^{-d_1-d_3}&&\bullet
1\ar@/^/[ll]^{rd+d_3}\ar@(dr,ur)_{d}}
$$


\subsection{Positivity and tightness}

\begin{prop}
The algebra $D(2B)^{r,c}$ is positively graded for every $c$ and
every $r$. The algebra $D(2B)^{r,c}$ is tightly graded if and only
if $c=0$ and $r=3$.
\end{prop}
\noindent {\bf Proof.} That $D(2B)^{r,c}$ is a positively graded
algebra follows easily from its relations. If ${\rm
deg}(\alpha)=2r$, ${\rm deg}(\beta)={\rm deg}(\gamma)=r$ and ${\rm
deg}(\eta)=4$, then the relations  are homogeneous.

If $D(2B)^{r,c}$ is tightly graded, then by Lemma \ref{lm: tight},
there are degree 1 elements of the form
$$
\begin{array}{c@{:=}l}
t_{\alpha}&\alpha+a_1\beta\gamma+a_2\alpha\beta\gamma,\\
t_{\beta}&\beta +b_1\alpha\beta,\\
t_{\gamma}&\gamma+b_2\gamma\alpha,\\
t_{\eta}&\eta+\sum_{i=2}^rd_i\eta^i,
\end{array}
$$
where $a_1,a_2,b_1,b_2,d_1,\dots,d_r$ are scalars.

It follows that $\alpha^2=t_{\alpha}^2$ is a homogeneous element
of degree 2, and that $\alpha\beta\gamma$ is a homogeneous element
of degree 3. If $c=1$, then this leads us to a contradiction. If
$c=0$, then from $\gamma\alpha\beta=t_{\gamma}t_{\alpha}t_{\beta}$
and $\eta^r=t_{\eta}^r$, we have that $r=3$. $\blacksquare$

\begin{prop}
The algebra $D(2A)^{r,0}$ is tightly graded for every $r$. The
algebra $D(2A)^{r,1}$ is positively graded if and only if $r\leq
2$. The algebras $D(2A)^{1,1}$ and $D(2A)^{2,1}$ are not tightly
graded.
\end{prop}
\noindent {\bf Proof.} If $r=0$, then it is obvious that if we put
the arrows of the quiver of $D(2A)^{r,0}$ in degree 1, then the
relations are homogeneous. Hence, $D(2A)^{r,0}$ is tightly graded.

If $c=1$ and  $r=1$, then if ${\rm deg}(\alpha)=2$, ${\rm
deg}(\beta)=1$ and ${\rm deg}(\gamma)=1$ we get a positive grading
on $D(2A)^{1,1}$. If $c=1$ and  $r=2$, then if ${\rm
deg}(\alpha)=2$, ${\rm deg}(\beta)=0$ and ${\rm deg}(\gamma)=0$,
we get a positive grading on $D(2A)^{2,1}$.

For $r>2$, if ${\rm deg}(\alpha)=r$, ${\rm deg}(\beta)=-(r-2)$ and
${\rm deg}(\gamma)=0$, we get a grading on $D(2A)^{r,1}$. The
graded quiver is given by
$$
\xymatrix{0\bullet\ar@(dl,ul)^{r}\ar@/^/[rr]^{2-r}&&\bullet
1\ar@/^/[ll]^{0}}
$$
This is not a positive grading. Also, this grading is not graded
Morita equivalent to the trivial grading on $D(2A)^{r,1}$. By
Proposition \ref{dvamodulajunik}, every other grading on
$D(2A)^{r,1}$ can be obtained from this grading by rescaling and
graded Morita equivalence. When we rescale a grading such that
there are homogeneous elements in both negative and positive
degrees, the resulting grading still has the same property. Let
$n_0$ and $n_1$ be integers and let ${\rm
Endgr}_{D(2A)^{r,1}}(P_0\langle n_0\rangle \oplus P_1\langle
n_1\rangle )^{op}$ be a graded algebra that is graded Morita
equivalent to the above graded algebra. By Proposition 9.1 in 
[\ref{GBTA}], the graded quiver of ${\rm
Endgr}_{D(2A)^{r,1}}(P_0\langle n_0\rangle\oplus P_1\langle
n_1\rangle)^{op}$ is given by
$$
\xymatrix{0\bullet\ar@(dl,ul)^{r}\ar@/^/[rr]^{(2-r)+n_0-n_1}&&\bullet
1\ar@/^/[ll]^{n_1-n_0}}
$$

If $(2-r)+n_0-n_1\geq 0$, then $n_1-n_0<0$. If $n_1-n_0\geq 0$,
then $(2-r)+n_0-n_1<0$. It follows that the resulting grading is
not positive. Hence, if $r>2$, then $D(2A)^{r,1}$ is not
positively graded.

To prove that $D(2A)^{2,1}$ is not tightly graded we start with
the grading on $D(2A)^{2,1}$ given by the graded quiver
$$
\xymatrix{0\bullet\ar@(dl,ul)^{1}\ar@/^/[rr]^{0}&&\bullet
1\ar@/^/[ll]^{0}}
$$
This grading is not graded Morita equivalent to the trivial
grading on $D(2A)^{2,1}$. As above, it follows easily that any
other grading that is graded Morita equivalent to this grading is
not positive. Hence, $D(2A)^{2,1}$ is not a tightly graded
algebra.

To prove that $D(2A)^{1,1}$ is not tightly graded we again use
Lemma \ref{lm: tight}. Assuming that $D(2A)^{1,1}$ is tightly
graded, we get that $\alpha^2$ is a homogeneous element of both
degree 2 and degree 3, which is impossible. $\blacksquare$

\section{One simple module}Any block with a dihedral defect group and one isomorphism class
of simple modules is Morita equivalent to some algebra from the
following family (cf.\ [\ref{ErK}] or [\ref{H96}]):

For a given integer $r\geq 1$, let $D:=D(1C)^r$ be the algebra
defined by the quiver and relations
\begin{multicols}{2}
$
\xymatrix{&&&&
\bullet\ar@(dl,ul)^{\alpha}\ar@(dr,ur)_{\beta}}
$

$\alpha^2=0=\beta^2$, $(\alpha\beta)^r=(\beta \alpha)^r.$
\end{multicols}

\subsubsection{Classification of gradings}

The relations of $D$ are homogeneous, regardless of the degrees of
$\alpha$ and $\beta$. It follows that for any pair of integers
$(a,b)$, we get a grading on $D$ by setting ${\rm deg}(\alpha)=a$
and ${\rm deg}(\beta)=b$. We denote this graded algebra by
$D^{a,b}$. When $a=b=1$ we get a tight grading on $D$. The graded
radical layers of the only projective indecomposable
$D^{a,b}$-module $_DD$ are
$$
\begin{array}{ccccc}
&&S&&\\
a&S&&S&b\\
a+b&S&&S&a+b\\
2a+b&S&&S&2b+a\\
&\vdots&&\vdots&\\
a+(r-1)(a+b)&S&&S&b+(r-1)(a+b)\\
&&S&&(a+b)^r\\
\end{array},
$$
where $S$ denotes the only simple $D$-module.

For a given integer $d$,  the graded
algebra ${\rm Endgr}_{D^{a,b}}(D\langle d \rangle)^{op}$ is graded
Morita equivalent to $D^{a,b}$ by Definition \ref{grMorDef}. But ${\rm
Endgr}_{D^{a,b}}(D\langle d \rangle)^{op}\cong D^{a,b}$, as graded
algebras. It follows that the only graded algebra which is graded
Morita equivalent to $D^{a,b}$ is $D^{a,b}$ itself. From this we
have the following proposition.

\begin{prop}
For any pair of integers $(a,b)$ there is a grading $D^{a,b}$ on $D$. For
different pairs of integers $(a,b)$ and $(c,d)$, the graded algebras
$D^{a,b}$ and $D^{c,d}$ are not graded Morita equivalent.
\end{prop}

It follows from this proposition that the maximal tori in ${\rm
Out }^K (D)$ are isomorphic to $\textbf{G}_m^l$, where $l>1$. If
it were that $l\leq 1 $, then we would have a unique grading up to
rescaling and graded Morita equivalence on $D$, which is not the
case.

If $\varphi$ is an arbitrary automorphism in ${\rm Out}^K(D)$,
then we can assume that
$$
\begin{array}{c@{=}l}
\varphi(e)&e,\\
\varphi(\alpha)&a_1\alpha+a_2\beta +a_3x,\\
\varphi(\beta)&b_1\alpha+b_2\beta +b_3y,
\end{array}
$$
where $a_i,b_i\in k$, and $x,y\in {\rm rad}^2D$. Since
$\varphi(\alpha^2)=\varphi(\beta^2)=0$, we have that $a_1a_2=0$
and $b_1b_2=0$. From $\varphi((\alpha\beta)^r)\neq 0$ and
$\varphi((\beta\alpha)^r)\neq 0$ it follows that either
$a_1\neq0\neq b_2$ and $a_2=b_1=0$, or $a_2\neq0 \neq b_1$ and
$a_1=b_2=0$. The action of $\varphi$ on ${\rm rad} D/{\rm rad}^2D$
is given by matrices of the form
$$\left(\begin{array}{cc}
a_1&0\\
0& b_2
\end{array}\right)\quad {\rm or}\quad \left(\begin{array}{cc}
0&b_1\\
a_2& 0
\end{array}\right).$$

It now follows easily (one can see this directly or by using
Remark 3.5 in [\ref{Rou}]) that the maximal tori in ${\rm Out}^K
(D)$ are isomorphic to the product of at most two copies of
$\textbf{G}_m$. Combining this conclusion with the above remarks
gives us that the maximal tori in ${\rm Out}^K(D)$ are isomorphic
to $\textbf{G}_m^2.$

\begin{prop}
The maximal tori in ${\rm Out}^K(D)$ are isomorphic to
$\textbf{G}_m^2.$ Up to graded Morita equivalence the gradings on
$D$ are parameterized by $\mathbb{Z}^2$ and are in one-to-one
correspondence with algebraic group homomorphisms from
$\textbf{G}_m$ to $\textbf{G}_m\times\textbf{G}_m$.
\end{prop}
\noindent {\bf Proof.} Follows from the above discussion and the
previous proposition. $\blacksquare$

\section{Summary of the results}

In the following table we summarize the results of this paper. The first three columns tell us respectively if there exists a non-trivial, a positive and a tight grading on a given block. The last column gives the isomorphism class of the maximal tori in the group of outer automorphisms of a given block. Derived
equivalence classes are separated by horizontal lines.
{\renewcommand{\arraystretch}{1.3}
$$
\begin{array}{c|ccccc}
{\rm Block}&{\rm Non}\!\!-\!\!{\rm trivial}&\phantom{aaa}{\rm Positive}\phantom{aaa}&\phantom{aaa}{\rm Tight }\phantom{aaa}&{\rm Maximal\,\, torus\,\,}&\\ \hline
A_r&{\rm Yes}&{\rm Yes}&{\rm Yes}&\textbf{G}_m\times \textbf{G}_m&\\
B_r&{\rm Yes}&{\rm Yes}&{\rm Only\,\, if}\, r=1&\textbf{G}_m\times \textbf{G}_m&\\
C_r&{\rm Yes}&{\rm Yes}&{\rm Only\,\, if}\, r=4&\textbf{G}_m\times \textbf{G}_m&\\ \hline
D(2A)^{r,0}&{\rm Yes}&{\rm Yes}&{\rm Yes}&\textbf{G}_m\times\textbf{G}_m&\\
D(2B)^{r,0}&{\rm Yes}&{\rm Yes}&{\rm Only\,\, if}\, r=3&\textbf{G}_m\times \textbf{G}_m&\\ \hline
D(2A)^{r,1}&{\rm Yes}&{\rm Only\,\, if}\, r\leq2&{\rm No}&\textbf{G}_m&\\
D(2B)^{r,1}&{\rm Yes}&{\rm Yes}&{\rm No}&\textbf{G}_m&\\ \hline
D(1C)^r&{\rm Yes}&{\rm Yes}&{\rm Yes}&\textbf{G}_m\times\textbf{G}_m&\\ \end{array}$$
}


\end{document}